\documentclass[12pt,leqno]{article}
\usepackage{amsmath,amssymb,amscd}
\title{On Rao's Theorems and the Lazarsfeld--Rao Property}
\author{Robin Hartshorne}
\date{}

\begin{document}

\maketitle

\begin{abstract}
Let $X$ be an integral projective scheme satisfying the condition $S_3$ of Serre
and $H^1({\mathcal O}_X(n)) = 0$ for all $n \in {\mathbb Z}$.  We generalize
Rao's theorem by showing that biliaison equivalence classes of codimension two
subschemes without embedded components are in one-to-one correspondence with 
pseudo-isomorphism classes of coherent sheaves on $X$ satisfying certain depth
conditions.

We give a new proof and generalization of Strano's strengthening of the
Lazarsfeld--Rao property, showing that if a codimension two subscheme is not
minimal in its biliaison class, then it admits a strictly descending elementary
biliaison.

For a three-dimensional arithmetically Gorenstein scheme $X$, we show that
biliaison equivalence classes of curves are in one-to-one correspondence with
triples $(M,P,\alpha)$, up to shift, where $M$ is the Rao module, $P$ is a
maximal Cohen--Macaulay module on the homogeneous coordinate ring of $X$, and
$\alpha: P^{\vee} \rightarrow M^* \rightarrow 0$ is a surjective map of the
duals.
\end{abstract}

\setcounter{section}{-1}
\section{Introduction}
\label{sec0}

The Lazarsfeld--Rao property was initially formulated for curves in ${\mathbb
P}^3$.  It says that each even liaison class of curves has minimal elements,
which form an irreducible family.  Any curve in the even liaison class that is
not minimal can be obtained from a minimal curve by a finite sequence of
ascending basic double links, followed by a deformation.  This result has been
generalized in several directions, and most recently Strano \cite{S} has shown
that if one uses elementary biliaisons, the deformation at the last step is not
necessary.  Note that a basic double link is essentially a union of the smaller
scheme with a complete intersection, while an elementary biliaison has linear
equivalence built into the definition.  So another way of phrasing Strano's
result is that the only deformations needed in the Lazarsfeld--Rao theorem are
linear equivalences on suitable hypersurfaces.

Our purpose in this paper is to give a uniform treatment of the theorems of Rao
and the Lazarsfeld--Rao property, in a fairly general context, and from a
slightly different point of view than the traditional one.  In particular, we
will use biliaisons only, and we deal with codimension~$2$ subschemes that are
equidimensional and without embedded components.

A good summary of the literature can be found in Chapter~$6$ of Migliore's book
\cite{M}.  Rao first gave theorems relating biliaison equivalence classes of
codimension~$2$ subschemes to stable equivalence classes of certain sheaves on
${\mathbb P}^n$ \cite{R1}, \cite{R2}.  The Lazarsfeld--Rao property was first
proved in a special case in \cite{LR}.  Then it was proved independently for
curves in ${\mathbb P}^3$ in \cite{MDP} and for codimension~$2$ locally
Cohen--Macaulay schemes in \cite{BBM}.  The ambient scheme was generalized from
${\mathbb P}^n$ to an arithmetically Gorenstein scheme in \cite{R2} and
\cite{BM}.  The condition that the subschemes be locally Cohen--Macaulay was
dropped in \cite{Na} and \cite{No}.  The necessity of the deformation at the
last step of the LR property was eliminated by Strano \cite{S}.

In this paper we retain all these generalizations.  Using ideas from the papers
mentioned above, and some technical results from the papers \cite{HMDP1},
\cite{HMDP2}, we can now prove the Lazarsfeld--Rao property for biliaison
equivalence classes of codimension~$2$ subschemes without embedded components
on an integral projective scheme $X$ satisfying condition $S_3$ of Serre and
$H_*^1({\mathcal O}_X) = 0$ $(2.4)$.  We also give a version of the theorem of
Rao for a three-dimensional arithmetically Gorenstein scheme $(3.2)$.  We have
made an effort in this paper to present each result with the minimal number of
hypotheses, so as to bring out more clearly which results depend on which
assumptions.

\section{Pseudo-isomorphism}
\label{sec1}

Here we borrow some ideas from \cite[Section~2]{HMDP1} and adapt them to our
situation.

\bigskip
\noindent
{\bf Hypotheses 1.1.}  Let $X$ be a projective scheme over an algebraically
closed field $k$, equidimensional of dimension $N \ge 2$.  We assume $X$
satisfies the condition $S_2$ of Serre.  Let ${\mathcal O}_X(1)$ be an ample
invertible sheaf.  For any coherent sheaf ${\mathcal F}$ we denote $\oplus_{n
\in {\mathbb Z}} H^i(X,{\mathcal F}(n))$ by $H_*^i({\mathcal F})$.  We assume
that
$H_*^1({\mathcal O}_X) = 0$.

\bigskip
\noindent
{\bf Definition 1.2.}  A coherent sheaf ${\mathcal L}$ on $X$ is {\em dissoci\'e}
if ${\mathcal L} \cong \oplus {\mathcal O}_X(n_i)$ for some $n_i \in {\mathbb
Z}$.

\bigskip
\noindent
{\bf Definition 1.3.}  An {\em elementary pseudo-isomorphism} is a surjective
map ${\mathcal E} \rightarrow {\mathcal E}' \rightarrow 0$ of coherent sheaves
whose kernel is dissoci\'e.  The equivalence relation generated by the
elementary pseudo-isomorphisms (and their inverses) is called {\em
pseudo-isomorphism} (psi for short).

\bigskip
\noindent
{\bf Lemma 1.4.} {\em A composition of elementary {\em psi} is an elementary
{\em psi}}.

\bigskip
\noindent
{\em Proof.}  Suppose given
\begin{eqnarray*}
& &0 \rightarrow {\mathcal L} \rightarrow {\mathcal E} \rightarrow {\mathcal E}'
\rightarrow 0 \\
& &0 \rightarrow {\mathcal L}' \rightarrow {\mathcal E}' \rightarrow {\mathcal
E}'' \rightarrow 0
\end{eqnarray*}
with ${\mathcal L}$ and ${\mathcal L}'$ dissoci\'e.  The kernel of the composed
map ${\mathcal E} \rightarrow {\mathcal E}''$ is an extension of ${\mathcal L}'$
by ${\mathcal L}$.  Because of the hypothesis $(1.1)$ that $H_*^1({\mathcal
O}_X) = 0$, this extension splits, so the kernel is ${\mathcal L} \oplus
{\mathcal L}'$, which is dissoci\'e.

\bigskip
\noindent
{\bf Lemma 1.5.}  {\em If ${\mathcal E}_1$ and ${\mathcal E}_2$ are {\em
psi}-equivalent, then there exists a coherent sheaf ${\mathcal F}$ and
elementary {\em psi}'s ${\mathcal F} \rightarrow {\mathcal E}_1$ and
${\mathcal F}
\rightarrow {\mathcal E}_2$.}

\bigskip
\noindent
{\em Proof.}  Any psi is a finite composition of elementary  psi's and their
inverses.  In view of $(1.4)$, and using induction on the length of a chain, it
is sufficient to show that if ${\mathcal E}_1 \rightarrow {\mathcal E}'$ and
${\mathcal E}_2 \rightarrow {\mathcal E}'$ are elementary psi's, then there
exists elementary psi's ${\mathcal F} \rightarrow {\mathcal E}_1$ and ${\mathcal
F} \rightarrow {\mathcal E}_2$ making a commutative diagram.  So let
\[
0 \rightarrow {\mathcal L}_1 \rightarrow {\mathcal E}_1 \rightarrow {\mathcal
E}' \rightarrow 0
\]
and
\[
0 \rightarrow {\mathcal L}_2 \rightarrow {\mathcal E}_2 \rightarrow {\mathcal
E}' \rightarrow 0.
\]
Take ${\mathcal F}$ to be the fibered sum $\ker({\mathcal E}_1 \oplus {\mathcal
E}_2 \rightarrow {\mathcal E}')$.  Then the kernel of the natural map ${\mathcal
F} \rightarrow {\mathcal E}_i$ is just ${\mathcal L}_{3-i}$ for $i = 1,2$.

\bigskip
\noindent
{\bf Definition 1.6.}  We say a coherent sheaf ${\mathcal E}$ on $X$ satisfies
the condition $T$ if the following hold

\begin{itemize}
\item[1)] ${\mathcal E}$ is locally free of constant rank in codimension $\le 1$.
\item[2)] ${\mathcal E}$ has depth $\ge 1$ in codimension $2$.
\item[3)] ${\mathcal E}$ has depth $\ge 2$ in codimension $\ge 3$.
\item[4)] There is a closed subset $Z \subseteq X$ of codimension $\ge 2$ such
that ${\mathcal E}$ is locally free on $X-Z$ and $\mbox{det } {\mathcal E}|_{X-Z}
\cong {\mathcal O}_{X-Z}(\ell)$ for some $\ell \in {\mathbb Z}$.  In this case,
we say ${\mathcal E}$ is {\em orientable}.
\end{itemize}

\bigskip
\noindent
{\bf Example 1.7.}  a) Let $V$ be a closed subscheme of $X$ of codimension $2$,
equidimensional, and having no embedded components.  Then the ideal sheaf
${\mathcal I}_V$ of $V$ satisfies the condition $T$.  To each such closed
subscheme $V$ then, we will associate the psi-equivalence class of coherent
sheaves satisfying $T$ containing the ideal sheaf ${\mathcal I}_V$.

b) Conversely, if ${\mathcal E}$ is a rank $1$ coherent sheaf satisfying $T$,
then ${\mathcal E}$ is isomorphic to ${\mathcal I}_V(n)$ for some subscheme $V$
of codimension $2$ without embedded points, and some $n \in Z$.  Indeed, let $Z$
be a closed set of codimension $\ge 2$ such that ${\mathcal E}$ is locally free
on $X-Z$.  Then since ${\mathcal E}$ is orientable, ${\mathcal E} \cong {\mathcal
O}_X(n)$ on $X-Z$ for some $n$.  Let $j: X-Z \rightarrow X$ be the inclusion. 
Then there is a natural map $\alpha: {\mathcal E} \rightarrow j_*j^*{\mathcal
E}$.  Since $X$ satisfies $S_2$, $j_*j^*{\mathcal E} \cong {\mathcal O}_X(n)$. 
Furthermore, the map $\alpha$ is injective because of the depth conditions on
${\mathcal E}$.  Hence ${\mathcal E} \cong {\mathcal I}_V(n)$ for some closed
subscheme $V$ of codimension $\ge 2$.  It may happen that ${\mathcal E} \cong
{\mathcal O}(n)$, in which case $V$ is empty.  If $V$ is nonempty, the depth
conditions on ${\mathcal E}$ now imply that $V$ is equidimensional of
codimension $2$, and has no embedded components.

c) If $V$ is a closed subscheme as in a) and b) above, and if $Y$ is a
hypersurface of $X$, linearly equivalent to $nH$, containing $V$, where $H$ is
the divisor class corresponding to ${\mathcal O}_X(1)$, then the exact sequence
$0
\rightarrow {\mathcal I}_Y \rightarrow {\mathcal I}_V \rightarrow {\mathcal
I}_{V,Y} \rightarrow 0$ shows that ${\mathcal I}_V \rightarrow {\mathcal
I}_{V,Y}$ is a psi, since ${\mathcal I}_Y \cong {\mathcal O}_X(-n)$.  Note,
however, that ${\mathcal I}_{V,Y}$ does not satisfy $T$, because it is a torsion
sheaf with support in codimension $1$ along $Y$.

\bigskip
\noindent
{\bf Lemma 1.8.}  a) {\em If ${\mathcal E} \rightarrow {\mathcal E}'$ is an
elementary {\em psi}, and ${\mathcal E}'$ satisfies $T$, then so does ${\mathcal
E}$.}

b) {\em if $X$ satisfies condition $S_3$ of Serre, and if ${\mathcal E}$ is a
sheaf satisfying $T$ of rank $\ge 2$, then there exists an elementary {\em psi}
${\mathcal E} \rightarrow {\mathcal E}'$ to a sheaf ${\mathcal E}'$ of rank one
less, also satisfying $T$.}

\bigskip
\noindent
{\em Proof.}  Part a) is obvious.  To prove b), let ${\mathcal E}(n)$ be a twist
of ${\mathcal E}$ that is generated by global sections.  I claim that if $s$ is
a sufficiently general section in $H^0({\mathcal E}(n))$, then the quotient
${\mathcal E}' = {\mathcal E}/{\mathcal O}(-n)$ defined by $s$ will be locally
free in codimension $1$.  At a point $x \in X$, the condition for ${\mathcal
E}'$ to be locally free is that $s(x) \in {\mathcal E}(n) \otimes k(x)$ be
non-zero.  If ${\mathcal E}$ is locally free of rank $r \ge 2$, since ${\mathcal
E}(n)$ is generated by global sections, the bad set of sections $s \in V =
H^0({\mathcal E}(n))$ will have codimension $r$.  It follows that for a general
$s \in V$, the bad set of $x \in X$ for which $s(x) = 0$ will have codimension
$\ge r \ge 2$.  Hence ${\mathcal E}'$ is locally free in codimension $1$.  The
rest of condition $T$ for ${\mathcal E}'$ follows from the $S_3$ condition on
$X$.

\bigskip
\noindent
{\bf Definition 1.9} \cite[2.6]{HMDP1}.  A coherent sheaf ${\mathcal E}$ $X$ is
{\em extraverti} if $H_*^1({\mathcal E}^{\vee}) = 0$ and ${\mathcal
E}xt^1({\mathcal E},{\mathcal O}_X) = 0$.

\bigskip
\noindent
{\bf Lemma 1.10.}  a) {\em If ${\mathcal E} \rightarrow {\mathcal E}'$ is an
elementary {\em psi}, and ${\mathcal E}'$ is extraverti, so is ${\mathcal E}$.}

b) {\em If ${\mathcal E}_1$ and ${\mathcal E}_2$ are extraverti sheaves that are
{\em psi} equivalent, then they are {\em stably equivalent}, namely, there exist
dissoci\'e sheaves ${\mathcal L}$ and ${\mathcal M}$ such that ${\mathcal E}_1
\oplus {\mathcal L} \cong {\mathcal E}_2 \oplus {\mathcal M}$.}

\bigskip
\noindent
{\em Proof.} a) If $0 \rightarrow {\mathcal L} \rightarrow {\mathcal E}
\rightarrow {\mathcal E}' \rightarrow 0$ with ${\mathcal L}$ dissoci\'e, then we
have
\[
0 \rightarrow {\mathcal E}'{}^{\vee} \rightarrow {\mathcal E}^{\vee} \rightarrow
{\mathcal L}^{\vee} \rightarrow {\mathcal E}xt^1({\mathcal E}',{\mathcal O}_X)
\rightarrow {\mathcal E}xt^1({\mathcal E},{\mathcal O}_X) \rightarrow 0.
\]
Assuming that ${\mathcal E}xt^1({\mathcal E}',{\mathcal O}_X) = 0$, we get
${\mathcal E}xt^1({\mathcal E},{\mathcal O}_X) = 0$ and an exact sequence $0
\rightarrow {\mathcal E}'{}^{\vee} \rightarrow {\mathcal E}^{\vee} \rightarrow
{\mathcal L}^{\vee} \rightarrow 0$.  Now the exact sequence of cohomology shows
$H_*^1({\mathcal E}^{\vee}) = 0$.

b) In view of part a) and $(1.5)$, it is sufficient to show that if $0
\rightarrow {\mathcal L} \rightarrow {\mathcal E} \rightarrow {\mathcal E}'
\rightarrow 0$ with ${\mathcal E}'$ extraverti, then the sequence splits, so
${\mathcal E} \cong {\mathcal E}' \oplus {\mathcal L}$.  For this we use the low
terms of the spectral sequence of Ext:
\[
0 \rightarrow H_*^1({\mathcal E}'{}^{\vee}) \rightarrow
\mbox{Ext}_*^1({\mathcal E}',{\mathcal O}_X) \rightarrow H_*^0({\mathcal
E}xt^1({\mathcal E}',{\mathcal O}_X)) \rightarrow H_*^2({\mathcal E}'{}^{\vee})
\rightarrow \dots
\]
The first and third terms being zero, so is the second, so the sequence splits.

\bigskip
\noindent
{\bf Proposition 1.11.} {\em Every {\em psi} equivalence class of sheaves
satisfying
$T$ contains an extraverti sheaf satisfying $T$.}

\bigskip
\noindent
{\em Proof.} Let ${\mathcal E}$ be a coherent sheaf satisfying $T$.  Let
\[
0 \rightarrow {\mathcal G} \rightarrow {\mathcal M} \rightarrow {\mathcal E}
\rightarrow 0
\]
be a resolution with ${\mathcal M}$ dissoci\'e.  From the hypothesis
$H^1_*({\mathcal O}_X) = 0$ we see that $\mbox{Ext}_*^1({\mathcal E},{\mathcal
O}_X)$ is a quotient of $H_*^0({\mathcal G}^{\vee})$, and hence is a finitely
generated graded $S = H_*^0({\mathcal O}_X)$-module.

Take a set of generators $\xi_i \in \mbox{Ext}^1({\mathcal E},{\mathcal
O}_X(-a_i))$ of this module and let
\[
0 \rightarrow {\mathcal L} \rightarrow {\mathcal F} \rightarrow {\mathcal E}
\rightarrow 0
\]
be the corresponding extension, with ${\mathcal L} = \oplus {\mathcal
O}_X(-a_i)$.  Then ${\mathcal F}$ is psi-equivalent to ${\mathcal E}$ by
construction, and satisfies $T$ by $(1.8)$.  We will show that ${\mathcal F}$ is
extraverti.

First we apply the functor $\mbox{Hom}(\cdot,{\mathcal O}_X(*))$.  This gives
\[
0 \rightarrow H_*^0({\mathcal E}^{\vee}) \rightarrow H_*^0({\mathcal F}^{\vee})
\rightarrow H_*^0({\mathcal L}^{\vee}) \stackrel{\alpha}{\rightarrow}
\mbox{Ext}_*^1({\mathcal E},{\mathcal O}_X) \rightarrow \mbox{Ext}_*^1({\mathcal
F},{\mathcal O}_X) \rightarrow 0.
\]
Furthermore, the map $\alpha$ sends the generators of $H_*^0({\mathcal
L}^{\vee})$ to the chosen generators $\xi_i$ of $\mbox{Ext}_*^1({\mathcal
E},{\mathcal O}_X)$.  Hence $\mbox{Ext}_*^1({\mathcal F},{\mathcal O}_X) = 0$,
and we conclude from the spectral sequence mentioned in the proof of $(1.10)$
above that $H_*^1({\mathcal F}^{\vee}) = 0$.

From that same spectral sequence, and noting that $H_*^2({\mathcal E}^{\vee})$ is
zero in large enough degrees, by Serre vanishing, we see that the map
\[
\mbox{Ext}_*^1({\mathcal E},{\mathcal O}_X) \rightarrow H_*^0({\mathcal
E}xt^1({\mathcal E},{\mathcal O}_X))
\]
is surjective in large enough degrees.  Hence the images of the generators
$\xi_i$ also generate the sheaf ${\mathcal E}xt^1({\mathcal E},{\mathcal O}_X)$. 
Then from applying ${\mathcal H}om(\cdot,{\mathcal O}_X)$ to the sequence
defining ${\mathcal F}$ above, we find also ${\mathcal E}xt^1({\mathcal
F},{\mathcal O}_X) = 0$, and ${\mathcal F}$ is extraverti.

\bigskip
Combining the results of this section, we have the following relationship between
codimension two subvarieties and coherent sheaves on $X$.

\bigskip
\noindent
{\bf Proposition 1.12.}  {\em Let $X$ be an equidimensional projective scheme of
dimension $\ge 2$ satisfying $S_2$ and $H_*^1({\mathcal O}_X) = 0$.}

a) {\em For any codimension two subscheme $V$ without embedded points, we
associate to it the {\em psi}-equivalence class of coherent sheaves satisfying $T$
that contains the sheaf ${\mathcal I}_V$.}

b) {\em If furthermore $X$ satisfies $S_3$, then every {\em psi}-equivalence class
of sheaves satisfying $T$ contains a sheaf of the form ${\mathcal I}_{V'}(n)$,
where
$V'$ is a codimension two subscheme without embedded points, for some $n \in
{\mathbb Z}$.}

c) {\em Every {\em psi}-equivalence class of sheaves satisfying $T$ contains
extraverti sheaves, unique up to stable equivalence.}

\bigskip
\noindent
{\em Proof.} a) is $(1.7)$a.

b) is $(1.8)$ together with $(1.7)$b.

c) is $(1.11)$ plus $(1.10)$.

\bigskip
\noindent
{\bf Remark 1.13.}  If $X$ satisfies in addition condition $G_2$, Gorenstein in
codimension $2$, then the condition ${\mathcal E}xt^1({\mathcal E},{\mathcal
O}_X) = 0$ implies by local duality depth ${\mathcal E} \ge 2$ at every point of
codimension $2$, so that then ${\mathcal E}$ will satisfy condition $S_2$.  In
the presence of $G_1$, this makes ${\mathcal E}$ reflexive \cite[1.9]{GD}. 
Thus, the extraverti sheaves of $(1.12)$c are reflexive, and we recover the
construction of Nollet
\cite{No} and Nagel
\cite{Na} by another route.

Furthermore, if $X$ is nonsingular, and we take $V$ to be locally
Cohen--Macaulay of codimension $2$, then depth ${\mathcal I}_V \ge r-1$ in
codimension $r$, and the same applies to ${\mathcal E}$.  Then the condition
${\mathcal E}xt^1({\mathcal E},{\mathcal O}) = 0$ makes depth ${\mathcal E} = r$
in codimension $r$, i.e., ${\mathcal E}$ is a Cohen--Macaulay sheaf.  Since $X$
is nonsingular, the extraverti sheaves associated to $V$ will be locally free. 
Then we recover the usual $N$-type resolution of ${\mathcal I}_V$.

\section{The Lazarsfeld--Rao property}
\label{sec2}

We preserve the hypotheses $(1.1)$ and in addition we assume $X$ satisfies
$S_3$.  Let $H$ be the divisor class corresponding to ${\mathcal O}_X(1)$, i.e.,
whose fractional ideal is isomorphic to ${\mathcal O}_X(-1)$.  If $Y \subseteq
X$ is an effective divisor linearly equivalent to $nH$ for some $n$, then $Y$
is a Cartier divisor on $X$, and so the scheme $Y$ satisfies $S_2$.  Thus we can
speak of generalized divisors on $Y$ \cite{GDB}.

\bigskip
\noindent
{\bf Definition 2.1.}  Let $V$ be a closed subscheme of pure codimension $2$ of
$X$, with no embedded components.  We say another such subscheme $V'$ is
obtained by an {\em elementary biliaison} of {\em height} $h$ from $V$, if there
exists an effective divisor $Y \sim nH$ on $X$ for some $n$, and a linear
equivalence $V' \sim V + hH$ on $Y$ for some $h \in {\mathbb Z}$.  (Here by
abuse of notation, $H$ denotes the divisor class of ${\mathcal O}_X(1)$ or
${\mathcal O}_Y(1)$ according to context.)  The equivalence relation generated
by elementary biliaisons will be called {\em biliaison}.  (To be precise, we
should call this notion {\em complete intersection biliaison} or CI-{\em
biliaison}, to distinguish it from the more general notion of {\em Gorenstein
biliaison} that uses ACM divisors $Y$ on $X$ \cite{GDB}.)

\bigskip
\noindent
{\bf Proposition 2.2.}  {\em If $V_1$ and $V_2$ are closed subschemes of
codimension $2$ of $X$ as above that are equivalent for biliaison, then there is
an integer $m \in {\mathbb Z}$ such that the ideal sheaves ${\mathcal I}_{V_1}$
and ${\mathcal I}_{V_2}(m)$ are equivalent for pseudo-isomorphism.}

\bigskip
\noindent
{\em Proof.}  Let $V' \sim V + hH$ as a divisor $Y \sim nH$ of $X$.  Then
${\mathcal I}_{V',Y} \cong {\mathcal I}_{V,Y}(-h)$.  On the other hand,
${\mathcal I}_{V',Y} \sim {\mathcal I}_{V',X}$ for psi and ${\mathcal I}_{V,Y}
\sim {\mathcal I}_{V,X}$ for psi by $(1.7)$.  Hence ${\mathcal I}_{V'} \sim
{\mathcal I}_V(-h)$ for psi.  Combining a sequence of elementary biliaisons
gives the result.

\bigskip
\noindent
{\bf Remark 2.3.}  The converse of this proposition is the theorem of Rao:  if
the ideal sheaves ${\mathcal I}_{V_1}$ and ${\mathcal I}_{V_2}(n)$ are
equivalent for psi, then $V_1$ and $V_2$ are equivalent for biliaison.  The
Lazarsfeld--Rao property gives a more detailed structure of the biliaison
equivalence class of codimension $2$ subschemes.  Both results are combined in
the following theorem.

\bigskip
\noindent
{\bf Theorem 2.4.}  {\em Let $X$ be a projective scheme over an algebraically
closed field $k$, equidimensional of dimension $N \ge 2$.  Let ${\mathcal O}_X(1)$
be a very ample invertible sheaf.  We assume that $X$ is integral, that it
satisfies condition $S_3$ of Serre, and that $H_*^1({\mathcal O}_X) = 0$.  We
consider closed subschemes $V$ of codimension $2$ equidimensional and without
embedded components.}

\begin{itemize}
\item[a)] {\em If ${\mathcal I}_{V_1}$ and ${\mathcal I}_{V_2}(n)$ are {\em
psi}-equivalent for some $n \in {\mathbb Z}$, then $V_1$ and $V_2$ are
equivalent for biliaison.  Hence the biliaison equivalence classes of subschemes
$V$ are in one-to-one correspondence with {\em psi}-equivalence classes (up to
twist) of coherent sheaves satisfying condition $T$ $(1.6)$.}

\item[b)] {\em If $V$ is a codimension $2$ subscheme whose degree is not minimal
in its biliaison equivalence class, then $V$ admits a strictly descending
biliaison (i.e., there exists a divisor $Y \sim nH$ on $X$ containing $V$, and a
subscheme $V' \sim V + hH$ on $Y$ with $h < 0$).}

\item[c)] {\em Any two subschemes $V,V'$ in the same biliaison class, both of
minimal degree, can be joined by a sequence of elementary biliaisons of height
$0$, i.e., linear equivalences on divisors $Y_i \sim n_iH$ on $X$.}
\end{itemize}

\bigskip
The proof will follow after some preliminary results.  The main idea is this: 
if ${\mathcal I}_{V_1}(a)$ and ${\mathcal I}_{V_2}(b)$ are equivalent for psi,
then by $(1.5)$ there exists a coherent sheaf ${\mathcal E}$ (satisfying $T$ by
$(1.8)$) and elementary psi's from ${\mathcal E}$ to ${\mathcal I}_{V_1}(a)$ and
${\mathcal I}_{V_2}(b)$.  Thus there are exact sequences
\[
0 \rightarrow \oplus_{i=1}^r {\mathcal O}(-a_i) \stackrel{\alpha}{\rightarrow}
{\mathcal E} \rightarrow {\mathcal I}_{V_1}(a) \rightarrow 0
\]
\[
0 \rightarrow \oplus_{i=1}^r {\mathcal O}(-b_i) \stackrel{\beta}{\rightarrow}
{\mathcal E} \rightarrow {\mathcal I}_{V_2}(b) \rightarrow 0
\]
where ${\mathcal E}$ has rank $r+1$.  The maps $\alpha,\beta$ are defined by
sections $s_i\in H^0({\mathcal E}(a_i))$ and $t_i \in H^0({\mathcal E}(b_i))$. 
The proof proceeds by comparison of the integers $a_i,b_i$, and a careful study
of the exact conditions for a section $s \in H^0({\mathcal E}(n))$ to have a
quotient ${\mathcal E}' = {\mathcal E}(n)/(s)$ that satisfies $T$.

First we consider a special case.

\bigskip
\noindent
{\bf Proposition 2.5.} {\em Let $X$ satisfy the hypotheses of $(2.4)$, let
${\mathcal E}$ be a rank $2$ coherent sheaf on $X$ satisfying $T$, and suppose
there are codimension $2$ subschemes $V_1,V_2$ and exact sequences}
\[
0 \rightarrow {\mathcal O}(-a_1) \stackrel{\alpha}{\rightarrow} {\mathcal E}
\rightarrow {\mathcal I}_{V_1}(a) \rightarrow 0
\]
\[
0 \rightarrow {\mathcal O}(-b_1) \stackrel{\beta}{\rightarrow} {\mathcal E}
\rightarrow {\mathcal I}_{V_2}(b) \rightarrow 0.
\]
{\em Then $V_2$ is obtained from $V_1$ by an elementary biliaison of height
$b-a$ on a suitable divisor $Y \sim nH$ on $X$.}

\bigskip
\noindent
{\em Proof.}  Consider the composed map $\gamma: {\mathcal O}(-b_1)
\stackrel{\beta}{\rightarrow} {\mathcal E} \rightarrow {\mathcal I}_{V_1}(a)$. 
If $\gamma$ is zero, then $\beta$ factors through ${\mathcal O}(-a_1)$.  We get
an injective map ${\mathcal O}(-b_1) \rightarrow {\mathcal O}(-a_1)$ whose
cokernel is contained in ${\mathcal I}_{V_2}(b)$, hence is zero, because
${\mathcal I}_{V_2}(b)$ is torsion-free.  Thus $a_1 = b_1$, so $\alpha =
\beta$, $a = b$, and ${\mathcal I}_{V_1} = {\mathcal I}_{V_2}$, so $V_1$ and
$V_2$ are equal.

If $\gamma$ is not zero, then since $X$ is integral (and here is exactly where
we use the hypothesis $X$ integral), it is injective.  Let the cokernel be
${\mathcal F}$:
\[
0 \rightarrow {\mathcal O}(-b_1) \stackrel{\gamma}{\rightarrow} {\mathcal
I}_{V_1}(a) \rightarrow {\mathcal F} \rightarrow 0.
\]
Let $Y$ be the divisor on $X$ defined by $\gamma$ followed by the inclusion of
${\mathcal I}_{V_1}(a)$ in ${\mathcal O}_X(a)$.  Then $Y \sim (a+b_1)H$ on $X$,
and ${\mathcal F} \cong {\mathcal I}_{V_1,Y}(a)$.  A diagram chase shows that
${\mathcal F}$ also fits into another exact sequence
\[
0 \rightarrow {\mathcal O}(-a_1) \rightarrow {\mathcal I}_{V_2}(b) \rightarrow
{\mathcal F} \rightarrow 0.
\]
Hence ${\mathcal F} \cong {\mathcal I}_{V_2,Y}(b)$ also.  Therefore ${\mathcal
I}_{V_1,Y}(a) \cong {\mathcal I}_{V_2,Y}(b)$ and there is a linear equivalence
$V_2 \sim V_1 + (b-a)H$ on $Y$.  Thus $V_2$ is obtained by an elementary
biliaison of height $b-a$ from $V_1$.

\bigskip
\noindent
{\bf Proposition 2.6.}  {\em Let $X$ be an integral projective scheme
satisfying condition $S_2$ of Serre.  Let ${\mathcal E}$ be a torsion-free
coherent sheaf, locally free in codimension $1$.  Let $W$ be a subvector space
of $H^0({\mathcal E})$, and let 
${\mathcal E}_0$ be the subsheaf generated by $W$.  Then the following conditions
are equivalent:}

\begin{itemize}
\item[(i)] {\em There exists a section $s \in W$ such that
${\mathcal E}' = {\mathcal E}/(s)$ is torsion-free and locally free in
codimension $1$}
\item[(ii)] (a) {\em for all $x \in X$ of codimension $1$, {\em rank} $({\mathcal
E}_0 \otimes k(x) \stackrel{\sigma_x}{\rightarrow} {\mathcal E} \otimes k(x))
\ge 1$, and}
\item[] (b) {\em either {\em rank} ${\mathcal E}_0 \ge 2$ or ${\mathcal E}_0
\cong {\mathcal O}_X$ and ${\mathcal E}/{\mathcal E}_0$ is torsion-free and
locally free in codimension $1$.}
\end{itemize}

\bigskip
\noindent
{\em Proof} (cf. \cite[1.4]{HMDP2}) (i) $\Rightarrow$ (ii).  Suppose given $s
\in W$ such that ${\mathcal E}' = {\mathcal E}/(s)$ is
torsion-free and locally free in codimension $1$.  Then there is a diagram
\[
\begin{CD}
0 @>>> {\mathcal O} @>s>> {\mathcal E} @>>> {\mathcal E}' @>>> 0 \\
@.        @VVV               @VVV               @VVV \\
0 @>>> {\mathcal E}_0 @>>> {\mathcal E} @>>> {\mathcal E}/{\mathcal E}_0 @>>>
0.
\end{CD}
\]
Let $x$ be a point of codimension $1$, and tensor with $k(x)$.  Since ${\mathcal
E}'$ is locally free at $x$, we have
\[
\begin{CD}
0 @>>> k(x) @>>> {\mathcal E} \otimes k(x) @>>> {\mathcal E}' \otimes k(x) @>>>
0 \\
@. @VVV   @VVV   @VVV \\
@. {\mathcal E}_0 \otimes k(x) @>{\sigma_x}>> {\mathcal E} \otimes k(x) @>>>
({\mathcal E}/{\mathcal E}_0) \otimes k(x) @>>> 0.
\end{CD}
\]
Thus it is clear that rank $\sigma_x \ge 1$.  If rank ${\mathcal E}_0 \ge 2$,
there is nothing more to prove.  If rank ${\mathcal E}_0 = 1$, then ${\mathcal
E}_0/{\mathcal O} \subseteq {\mathcal E}'$, and since ${\mathcal E}'$ is
torsion-free, we must have ${\mathcal E}_0 \cong {\mathcal O}$ and ${\mathcal
E}/{\mathcal E}_0 \cong {\mathcal E}'$ torsion-free and locally free in
codimension $1$, as required.

(ii) $\Rightarrow$ (i).  Assuming condition (ii), we will show that a general
element $s \in W$ makes ${\mathcal E}' = {\mathcal E}/(s)$
locally free in codimension $1$.  Then the condition ${\mathcal E}$ torsion-free
and $X$ satisfies $S_2$ will imply ${\mathcal E}'$ also torsion-free.

Since ${\mathcal E}$ is locally free in codimension $1$, the condition that
${\mathcal E}'$ be locally free at a point $x$ of codimension $1$ is precisely
that $s(x) \in {\mathcal E} \otimes k(x)$ be non-zero.  Consider the bad locus $B
= \{(s,x) \mid s(x) = 0$ in ${\mathcal E}
\otimes k(x)\}$ inside $W \times X$.  From condition (ii)~a) it is clear that
${\mathcal E}_0 \ne 0$.  If rank ${\mathcal E}_0 = 1$, then we know that
${\mathcal E}_0 \cong {\mathcal O}_X$, and ${\mathcal E}/{\mathcal E}_0$ locally
free in codimension $1$, so there is nothing to prove.  So we may assume rank
${\mathcal E}_0 \ge 2$.

Since $X$ is integral, there is a nonempty open set $U \subseteq X$ over which
${\mathcal E}/{\mathcal E}_0$ is locally free.  For any $x \in U$, rank
$\sigma_x \ge 2$ and so the fiber $B_x$ of $B$ over $x$ has codimension $\ge 2$
in $W$.  Therefore $B_U = B \cap (W \times U)$ has dimension $\le w - 2 + N$,
where $w = \dim W$ and $N = \dim X$.  Let $D = X-U$.  At points $x \in D$ of
codimension $1$ in $X$, rank $\sigma_x \ge 1$.  Hence rank $\sigma_x \ge 1$ over
an open set $U' \le D$, and $\dim B_{U'} \le w - 1 + N - 1$.  Finally, there may
be some subset $Z \subseteq D_1$ of codimension $\ge 2$ in $X$, for which rank
$\sigma_x = 0$.  Then $\dim B_Z \le w + N - 2$.  All in all, we find $\dim B \le
w + N - 2$.  Looking at the projection of $B$ to $W$, we see that for a general
$s \in W$, the fiber $B_s \subseteq X$ of $B$ over $s$ will have codimension
$\ge 2$.  Therefore ${\mathcal E}' = {\mathcal E}/(s)$ will be locally free in
codimension $1$, as required.

\bigskip
\noindent
{\bf Remark 2.7.} If $X$ satisfies the hypotheses $(1.1)$ and also condition
$S_3$ of Serre, and if the sheaf ${\mathcal E}$ of $(2.6)$ satisfies $T$, then
the quotients ${\mathcal E}' = {\mathcal E}/(s)$ of $(2.6)$ will also satisfy
$T$.  Indeed the extra conditions 3), 4) of $(1.6)$ follow from the same
conditions for ${\mathcal E}$ plus the condition $S_3$ on $X$.

\bigskip
\noindent
{\bf Proposition 2.8.}  {\em With the hypotheses of $(2.4)$, suppose that
${\mathcal E}$ is a coherent sheaf satisfying $T$, and that there are
codimension $2$ subschemes $V$ and $W$ and exact sequences}
\[
0 \rightarrow \oplus_{i=1}^r {\mathcal O}(-a_i) \stackrel{\alpha}{\rightarrow}
{\mathcal E} \rightarrow {\mathcal I}_V(a) \rightarrow 0 
\]
\[
0 \rightarrow \oplus_{i=1}^r {\mathcal O}(-b_i) \stackrel{\beta}{\rightarrow}
{\mathcal E} \rightarrow {\mathcal I}_W(b) \rightarrow 0.
\]
{\em Assume that $a_1 \le a_2 \le \dots \le a_r$ and $b_1 \le b_2 \le \dots \le
b_r$ and that $a_i = b_i$ for $i = 1,\dots,k-1$, and $a_k < b_k$ for some $k$. 
Then $W$ admits a strictly descending biliaison.}

\bigskip
\noindent
{\em Proof.} Let $s_i \in H^0({\mathcal E}(a_i))$ and $t_i \in H^0({\mathcal
E}(b_i))$ be the sections defining the maps $\alpha$ and $\beta$.  The idea is to
show that $t_1,\dots,t_{k-1},s_k,t_{k+1},\dots,t_r$ define another subscheme
$W'$ by an exact sequence
\[
0 \rightarrow \oplus_{i\ne k} {\mathcal O}(-b_i) \oplus {\mathcal O}(-a_k)
\rightarrow {\mathcal E} \rightarrow {\mathcal I}_{W'}(b') \rightarrow 0.
\]
Then considering determinants of the locally free sheaves on an open subset of
$X$ we find $b' = b - b_k + a_k < b$.  If we let ${\mathcal E}' = {\mathcal
E}/(t_1,\dots,t_{k-1},t_{k+1},\dots,t_r)$, then ${\mathcal E}'$ is a rank $2$
coherent sheaf with sections $t_k,s_k$ defining $W$ and $W'$, respectively. 
According to $(2.5)$, $W'$ is then obtained from $W$ by an elementary biliaison
of height $b' - b < 0$, and we get the desired result.  This idea may not work
with the original $s_i,t_i$, so we must modify them a little, without changing
$V$ and $W$.

\bigskip
\noindent
{\bf Step 1.} We will show that for sufficiently general $s'_k \in H^0({\mathcal
E}(a_k))$, the quotient sheaf ${\mathcal E}/(t_1,\dots,t_{k-1},s'_k)$ satisfies
$T$ of rank $r+1-k$.  To show this, consider the sheaf ${\mathcal F} = {\mathcal
E}(a_k)/(t_1,\dots,t_{k-1})$.  This has ${\mathcal I}_W(b+a_k)$ as a quotient,
hence satisfies $T$ by $(1.8)$.  We will apply $(2.6)$ to ${\mathcal F}$ to show
that it has a section ${\bar s}'_k$ (which lifts to a section $s'_k \in
H^0({\mathcal E}(a_k))$) with quotient satisfying $T$.

Let ${\mathcal E}_0$ be the subsheaf of ${\mathcal E}(a_k)$ generated by global
sections, and let ${\mathcal F}_0$ be the subsheaf of ${\mathcal F}$ generated
by global sections.  Then we have
\[
\begin{CD}
0  @>>>  \oplus_{i=1}^{k-1} {\mathcal O}(a_k-b_i)  @>{(t_i)}>>  {\mathcal E}_0 
@>>>  {\mathcal F}_0  @>>>  0 \\
@. \| @.  @VVV     @VVV \\
0  @>>>  \oplus_{i=1}^{k-1} {\mathcal O}(a_k-b_i) @>{(t_i)}>> {\mathcal E}(a_k)
@>>> {\mathcal F} @>>> 0.
\end{CD}
\]
On the other hand, we know that ${\mathcal E}/(s_1,\dots,s_k)$ satisfies $T$ and
is of rank $r+1-k$, so for every point $x \in X$ of codimension $1$ we must have
rank $(\sigma_x({\mathcal E}(a_k)): {\mathcal E}_0 \otimes k(x) \rightarrow
{\mathcal E}(a_k)\otimes k(x)) \ge k$.  Tensoring the diagram above with $k(x)$,
it follows that rank $(\sigma_x({\mathcal F}): {\mathcal F}_0 \otimes k(x)
\rightarrow {\mathcal F} \otimes k(x)) \ge 1$.  Furthermore, if rank ${\mathcal
F}_0 = 1$, then rank ${\mathcal E}_0 = k$, so ${\mathcal E}_0$ must be equal to
the subsheaf generated by $s_1,\dots,s_k$, which is $\oplus_{i=1}^k {\mathcal
O}(a_k-a_i)$.  In this case ${\mathcal E}_0/(t_1,\dots,t_{k-1}) \cong {\mathcal
E}_0/(s_1,\dots,s_{k-1}) \cong {\mathcal O}_X$ and ${\mathcal F}/{\mathcal F}_0
= {\mathcal E}(a_k)/{\mathcal E}_0$ satisfies $T$.

So the conditions of $(2.6)$ are satisfied (taking the whole space $H^0({\mathcal
E}(a_k))$), and there is a section
$s'_k
\in {\mathcal E}(a_k)$ so that ${\mathcal E}(a_k)/(t_1,\dots,t_{k-1},s'_k) =
{\mathcal F}/({\bar s}'_k)$ satisfies $T$.  Note also that the image of $s'_k$
in ${\mathcal I}_W(b+a_k)$ is nonzero.  For if it were zero, then $s'_k \in
H^0\left( \oplus_{i=1}^r (a_k-b_i)\right)$.  Since $a_k < b_k$, in fact $s'_k
\in H^0\left( \oplus_{i=1}^{k-1} {\mathcal O}(a_k-b_i)\right)$, which is
impossible since ${\bar s}'_k \ne 0$ in $H^0({\mathcal F})$.

\bigskip
\noindent
{\bf Step 2} (cf.\ proof of \cite[2.1]{BBM}).  Now we will show that for
suitable choice of $t'_i = t_i + f_it_k$, $i = k+1,\dots,r$, where the $f_i$ are
elements of $H^0({\mathcal O}_X(b_i-b_k))$, the sheaf ${\mathcal E}'' =
{\mathcal E}/(t_1,\dots,t_{k-1},s'_k,t'_{k-1},\dots,t'_r)$ will satisfy $T$ of
rank $1$, hence be of the form ${\mathcal I}_{W'}(b')$ for some codimension $2$
scheme $W'$.

To do this, consider the sheaf ${\mathcal G} = {\mathcal E}/(t_1,\dots,t_r,s'_k)
= {\mathcal E}''/(t_k)$.  Since ${\mathcal G}$ is the quotient of ${\mathcal
I}_W(b)$ by the image of $s'_k$, it is a torsion sheaf supported along a divisor
$Y$.  Thus ${\mathcal E}''$ is locally free in codimension $1$ except possibly
along $Y$.  For each generic point $y$ of an irreducible component of $Y$, we
have $t_1(y),\dots,t_{k-1}(y),s'_k(y)$ linearly independent in ${\mathcal E}
\otimes k(y)$, and similarly for $t_1(y),\dots,t_r(y)$.  Hence for general forms
$f_i \in H^0({\mathcal O}_X(b_i-b_k))$, letting $t'_i = t_i + f_it_k$, $i =
k+1,\dots,r$, we will have
$t_1(y),\dots,t_{k-1}(y),s'_k(y),t'_{k+1}(y),\dots,t'_r(y)$ linearly
independent, and so ${\mathcal E}''$ will be locally free at $y$.  We can do
this simultaneously for the finite number of generic points $y$ of $Y$.  Thus
${\mathcal E}''$ will be locally free of rank $1$ in codimension $1$, hence will
satisfy $T$ and be of the form ${\mathcal I}_{W'}(b')$.

\bigskip
\noindent
{\bf Step 3.}  Now take ${\mathcal E}' = {\mathcal
E}/(t_1,\dots,t_{k-1},t'_{k+1},\dots,t'_r)$.  Then ${\mathcal E}'$ satisfies $T$
of rank $2$, and ${\mathcal E}'/(t_k) = {\mathcal I}_W(b)$, and ${\mathcal
E}'/(s'_k) = {\mathcal I}_{W'}(b')$.  Then, as explained above, using $(2.5)$,
$W'$ is obtained by a strictly descending elementary biliaison from $W$.

\bigskip
\noindent
{\bf Proposition 2.9.}  {\em Suppose given $V,W$ and exact sequences as in
$(2.8)$ with $a_i = b_i$ for each $i$.  Then $W$ is obtained from $V$ by a
finite number of elementary biliaisons of height zero.}

\bigskip
\noindent
{\em Proof.} By induction on rank ${\mathcal E} = r+1$.  Note first that by
computing determinants, since $a_i = b_i$ for all $i$, we find $a=b$.  If rank
${\mathcal E} = 2$, this is just $(2.5)$.  So suppose rank ${\mathcal E} \ge 3$,
i.e., $r \ge 2$.  As in the proof of $(2.8)$, let $s_i$ and $t_i$ define the
maps $\alpha$ and $\beta$.

Applying $(2.6)$ to the sheaf ${\mathcal F} = {\mathcal E}/(s_1,\dots,s_{r-1})$,
since ${\mathcal E}(a_r)$ has one section $s_r$ making a quotient satisfying
$T$, it follows that if we take a general section $s' \in H^0({\mathcal
E}(a_r))$, then ${\mathcal E}' = {\mathcal E}/(s_1,\dots,s_{r-1},s')$ will be
equal to ${\mathcal I}_{V'}(a)$ for some codimension $2$ subscheme $V'$.  Then
using $(2.5)$ applied to ${\mathcal F}$, we find $V$ and $V'$ are related by an
elementary biliaison of height zero.

Doing the same with ${\mathcal G} = {\mathcal E}/(t_1,\dots,t_{r-1})$, we find
that for sufficiently general $s' \in H^0({\mathcal E}(a_r))$, which we can take
to be the same as the $s'$ above (!), the quotient ${\mathcal
E}/(t_1,\dots,t_{r-1},s')$ will define a subscheme $W'$, related by one
elementary biliaison of height zero to $W$.

Now $V'$ and $W'$ both have resolutions of the above form using the sheaf
${\mathcal E}' = {\mathcal E}/(s')$, of rank one less, so by the induction
hypothesis, $V'$ and $W'$ are related by a finite sequence of elementary
biliaisons of height zero, and we are done.

\bigskip
\noindent
{\em Proof of $(2.4)$.}  Suppose given subschemes $V$ and $W$ such that
${\mathcal I}_V(a)$ and ${\mathcal I}_W(b)$ are psi-equivalent.  Then by $(1.5)$
there is a coherent sheaf ${\mathcal E}$ satisfying $T$ and exact sequences as
in the statement of $(2.8)$.  If $a_i$ is not equal to $b_i$ for all $i$, then
by $(2.8)$, one of the two subschemes admits a descending biliaison to a
subscheme of lower degree.  Since the degree of a subscheme is always
nonnegative, we can proceed inductively, and after a finite number of steps we
arrive at a situation where $a_i = b_i$ for all $i$.  Then $(2.9)$ applies to
show they are equivalent by a finite  number of biliaisons of height zero.

This proves a) of the theorem.  For b), given $V$ not of minimal degree, take
any $W$ of minimal degree in the same biliaison class.  Then in $(2.8)$ it must
be $V$ that admits the descending biliaison.

As for c), if $V$ and $W$ are both of minimal degree, neither one can admit a
descending biliaison, so the $a_i$ must be equal to the $b_i$, and $(2.9)$
applies.

\bigskip
\noindent
{\bf Remark 2.10.}  There is one special case of the theorem that merits special
attention.  It might happen in the course of the proof that one of the
subschemes obtained by a descending biliaison is empty.  The theorem and its
proof still hold, provided that we allow the empty scheme.  The psi class of the
corresponding coherent sheaves is the class containing the dissoci\'e sheaves. 
This corresponds to one biliaison equivalence class, the one containing complete
intersections of hypersurfaces $Y \sim nH$ and $Y' \sim n'H$ in $X$.  The
schemes of minimal degree in this class are the empty scheme.

This may be regarded as unsatisfactory, so we prove separately that for the
biliaison class corresponding to the dissoci\'e sheaves, the results of the
theorem hold also if we restrict our attention only to the nonempty schemes.

We have only to consider the case where in the proof of $(2.8)$ we might obtain
an empty scheme $W'$.  In that case, $W \sim mH$ on a hypersurface $Y \sim nH$
in $X$.  Since $H^1({\mathcal O}_X(m-n)) = 0$, it follows that $H^0({\mathcal
O}_X(m)) \rightarrow H^0({\mathcal O}_Y(m))$ is surjective, so $W$ is a complete
intersection of $Y$ with $Y' \sim mH$.  If $m \ge 2$, we can make a descending
biliaison to $W' \sim H$ on $Y$.  If $m=1$ and $n\ge 2$, we regard $W$ as a
divisor on $Y' \sim H$ on $X$ and again make a descending biliaison to $W' \sim
H$ on $Y'$.  If $m=n=1$, I claim $W$ is of minimal degree among the nonempty
schemes in this biliaison class.  For if $V$ was something of lower degree,
according to the theorem there would be a sequence $V = V_1,V_2,\dots,V_r$ of
descending biliaisons, with $V_r$ the empty scheme.  In that case, as we have
seen, $V_{r-1}$ would be a complete intersection of degree $\ge W$, which is
impossible.  The same argument shows that any nonempty scheme of minimal degree
is of the same form as $W$.  Now $W$ has a resolution of the form $0 \rightarrow
{\mathcal O}(-2) \rightarrow {\mathcal O}(-1) \oplus {\mathcal O}(-1)
\rightarrow {\mathcal I}_W \rightarrow 0$.  Given two of these, they have the
same ${\mathcal E}$ and the same $a_i$, so by $(2.9)$ they are joined by a
sequence of biliaisons of height zero.

\bigskip
\noindent
{\bf Remark 2.11.}  In case $X = {\mathbb P}^n$, the codimension two schemes in
the biliaison class of the empty scheme are the ACM schemes.  So we have just
shown $(2.10)$ that any ACM codimension $2$ scheme admits descending elementary
biliaisons to a minimal one, which is just a linear variety ${\mathbb P}^{n-2}$
in ${\mathbb P}^n$.

\bigskip
\noindent
{\bf Remark 2.12.}  A consequence of the theorem is that a variety $V$ that is
not of minimal degree can be reached by a sequence of strictly ascending
elementary biliaisons from one of minimal degree.  For curves in ${\mathbb
P}^3$, this is Strano's theorem \cite{S}, answering a question of
\cite[p.~93]{MDP}, showing that the deformation in earlier proofs of the
Lazarsfeld--Rao property is not necessary.

\bigskip
\noindent
{\bf Remark 2.13.}  If we wish to prove only $(2.4a)$, namely that psi
equivalence of ideal sheaves implies biliaison equivalence of subschemes, we can
remove the hypothesis that $X$ is integral.  We prove this by the original method
of Rao.  First of all, in $(2.5)$ instead of comparing $\alpha$ to $\beta$, we
compare them both to $\gamma: {\mathcal O}(-c) \rightarrow {\mathcal E}$, where
$c$ is chosen sufficiently large that ${\mathcal E}(c)$ is generated by global
sections.  Then the induced maps ${\mathcal O}(-c) \rightarrow {\mathcal
I}_{V_1}(a)$ and ${\mathcal O}(-c) \rightarrow {\mathcal I}_{V_2}(b)$ will both
be injective, and we find that $V_1$ and $V_2$ are now related by two biliaisons
of large heights.

For sheaves of arbitrary rank, we follow the plan of proof of $(2.9)$, but take
the section $s'$ of that proof to be a general section of ${\mathcal E}(n)$ for
$n \gg 0$.  Then we can make $s'$ independent of the other sections and proceed
as in that proof, except that the biliaisons are now up and down of large
heights.

\bigskip
\noindent
{\bf Corollary 2.14.} a) {\em With the hypotheses of $(2.4)$, the biliaison
equivalence classes of codimension two subschemes are in one-to-one
correspondence with stable equivalence classes of extraverti sheaves satisfying
$T$ (up to twist).}

b) {\em If furthermore $X$ also satisfies $G_2$, then the biliaison classes of
codimension two subschemes are in one-to-one correspondence with stable
equivalence classes of reflexive extraverti sheaves satisfying $T$ (up to twist).}

\bigskip
\noindent
{\em Proof.} This is simply a restatement of $(2.4)$, using $(1.12)$ and
$(1.13)$.  Part b) is the theorem of Nollet \cite[2.12]{No} in the case $X =
{\mathbb P}^n$, and of Nagel \cite[6.4]{Na} in the case $X$ is an integral
arithmetically Gorenstein scheme.  Nagel also proved \cite[7.3]{Na} the weaker
form of the Lazarsfeld--Rao property (2.4b) allowing a deformation before making
a descending biliaison.

\section{Cohomological characterization of biliaison classes}
\label{sec3}

For curves in ${\mathbb P}^3$, the original theorem of Rao \cite{R1} says that
two curves $C_1,C_2$ are in the same biliaison equivalence class if and only if
their Rao modules $M_1,M_2$ are isomorphic up to twist.  The {\em Rao module} of
a curve $C$ is $M = H_*^1({\mathcal I}_C)$.  It was proved by first relating
biliaison classes of curves to stable equivalence classes of locally free
sheaves, and then using a theorem of Horrocks \cite{H}.  The theorem of Horrocks
is difficult to generalize to higher dimensions, but in this section we give a
version for an arithmetically Gorenstein scheme $X$ of dimension $3$.

\bigskip
\noindent
{\bf Proposition 3.1.} {\em Let $X$ be a $3$-dimensional projective {\em ACM}
scheme, and let $S = H_*^0({\mathcal O}_X)$.}

\begin{itemize}
\item[a)] {\em Each {\em psi} equivalence class of coherent sheaves ${\mathcal
E}$ satisfying $T$ determines a finite length graded $S$-module $M =
H_*^1({\mathcal E})$.}

\item[b)] {\em Conversely, for each finite length graded $S$-module $M$, there
is a locally free extraverti coherent sheaf ${\mathcal E}$ satisfying $T$, such
that $H_*^1({\mathcal E}) \cong M$.}
\end{itemize}

\bigskip
\noindent
{\em Proof.} For part a) we first show that if $0 \rightarrow {\mathcal L}
\rightarrow {\mathcal E} \rightarrow {\mathcal E}' \rightarrow 0$ is an
elementary psi, then $H_*^1({\mathcal E}) \cong H_*^1({\mathcal E}')$.  This is
because $X$ being ACM implies $H_*^i({\mathcal O}_X) = 0$ for $i = 1,2$. 
Serre's theorem tells us that $H^1({\mathcal E}(n)) = 0$ for $n \gg 0$.  To
show that $M$ is of finite length, we must show also $H^1({\mathcal E}(-n)) = 0$
for $n \gg 0$.  Let $\omega$ be the dualizing sheaf on $X$.  Then $H^1({\mathcal
E}(m))$ is dual to $\mbox{Ext}_X^2({\mathcal E},\omega(n))$.  For $n \gg 0$,
the sheaves ${\mathcal E}xt^i({\mathcal E},\omega(n))$ have no higher
cohomology, so the spectral sequence of local and global Ext gives
$\mbox{Ext}^2({\mathcal E},\omega(n)) = H^0({\mathcal E}xt^2({\mathcal
E},\omega(n)))$.  Now we use local duality at a closed point $x \in X$.  There
${\mathcal E}xt^2({\mathcal E},\omega)_x$ is dual to $H_x^1({\mathcal E})$, and
this is zero because depth ${\mathcal E}_x \ge 2$ by condition $T$.

For part b), let $M$ be a finite length graded $S$-module.  Take a resolution
\[
0 \rightarrow E \rightarrow L_1 \rightarrow L_0 \rightarrow M \rightarrow 0
\]
of graded $S$-modules, with $L_i$ free, and $E$ the second syzygy.  Passing to
associated sheaves we get
\[
0 \rightarrow {\mathcal E} \rightarrow {\mathcal L}_1 \rightarrow {\mathcal L}_0
\rightarrow 0.
\]
This shows that ${\mathcal E}$ is locally free and $H_*^1({\mathcal E}) \cong
M$.  Also $\mbox{det } {\mathcal E} = (\mbox{det } {\mathcal L}_1)(\mbox{det }
{\mathcal L}_0)^{-1}$, so ${\mathcal E}$ satisfies $T$.  Taking duals, we find
\[
0 \rightarrow {\mathcal L}_0^{\vee} \rightarrow {\mathcal L}_1^{\vee}
\rightarrow {\mathcal E}^{\vee} \rightarrow 0
\]
and ${\mathcal E}xt^1({\mathcal E},{\mathcal O}) = 0$.  From this, again using
$X$ is ACM, we find $H_*^1({\mathcal E}^{\vee}) = 0$, so ${\mathcal E}$ is
extraverti.

\bigskip
\noindent
{\bf Theorem 3.2.} {\em Now assume that $X$ is a $3$-dimensional arithmetically
Gorenstein scheme, i.e., $X$ is {\em ACM} and $S = H_*^0({\mathcal O}_X)$ is a
Gorenstein ring.}

a) {\em To each extraverti coherent sheaf ${\mathcal E}$ satisfying $T$
we will associate a maximal Cohen--Macaulay module $P$ over the ring $S$, defined
up to  stable equivalence, and a map $\alpha: P^{\vee} \rightarrow M^*
\rightarrow 0$, where
$M^*$ is the dual of the finite length module $M = H_*^1({\mathcal E})$.}

b) {\em If
${\mathcal E}_1$ and ${\mathcal E}_2$ are two such sheaves with the same
associated modules $M_1 \cong M_2$, and if the associated maximal
Cohen--Macaulay modules $P_1,P_2$ are stably equivalent and the maps $\alpha_1:
P_1^{\vee} \rightarrow M^*$ and $\alpha_2: P_2^{\vee} \rightarrow M^*$ are
compatible with the stable equivalence, then ${\mathcal E}_1$ and ${\mathcal
E}_2$ are stably equivalent.}

c) {\em Given a finite-length graded $S$-module $M$, a graded maximal
Cohen--Macaulay module $P$ whose associated sheaf ${\mathcal P}$ is orientable,
and given a map $\alpha: P^{\vee} \rightarrow M^* \rightarrow 0$, there exists an
extraverti sheaf ${\mathcal E}$ satisfying $T$ that gives rise to this triple, as
above.}

\bigskip
\noindent
{\em Proof.} a) Given ${\mathcal E}$ extraverti satisfying $T$, let $E =
H_*^0({\mathcal E})$, and take a resolution of graded $S$-modules
\[
0 \rightarrow P \rightarrow L_1 \rightarrow L_0 \rightarrow E \rightarrow 0
\]
where $L_0,L_1$ are free and $P$ is the kernel.  Since $E$ has depth $\ge 2$ at
the irrelevant prime ${\mathfrak m}$ of $S$, we see that $P$ has depth $4$, so
$P$ is a maximal Cohen--Macaulay module.  Taking duals we obtain an exact
sequence
\[
0 \rightarrow E^{\vee} \rightarrow L_0^{\vee} \rightarrow L_1^{\vee} \rightarrow
P^{\vee} \stackrel{\alpha}{\rightarrow} M^* \rightarrow 0.
\]
To see this, first consider $\mbox{Ext}_S^1(E,S)$.  This is dual on $S$ to
$H_{\mathfrak m}^3(E)$, which is isomorphic to $H_*^2({\mathcal E})$ on $X$. 
This in turn is dual to $\mbox{Ext}_{X,*}^1({\mathcal E},{\mathcal O})$
which is zero because ${\mathcal E}$ is extraverti, using the exact sequence of
low degree terms of the spectral sequence mentioned in the proof of $(1.10)$. 
Hence
$\mbox{Ext}_S^1(E,S) = 0$.  Secondly, $\mbox{Ext}_S^2(E,S)$ is dual to
$H_{\mathfrak m}^2(E)$, which is isomorphic to $H_*^1({\mathcal E})$, which is
$M$.

In this way we obtain the maximal CM module $P$, determined up to stable
equivalence on $S$, and the map $\alpha: P^{\vee} \rightarrow M^* \rightarrow 0$.

b) Now suppose ${\mathcal E}_1$ and ${\mathcal E}_2$ are two such sheaves, with
isomorphic associated modules $M_1 \cong M_2 = M$ and stably equivalent maximal
CM modules $P_1$ and $P_2$, and compatible maps $\alpha_1,\alpha_2$.  Then
$E_1^{\vee}$ and $E_2^{\vee}$ both occur in resolutions of $\alpha: P^{\vee}
\rightarrow M^* \rightarrow 0$ as above, and hence $E_1^{\vee}$ and $E_2^{\vee}$
are stably equivalent graded $S$-modules.  It follows that ${\mathcal E}_1$ and
${\mathcal E}_2$ are stably equivalent sheaves.

c) Given $\alpha: P^{\vee} \rightarrow M^* \rightarrow 0$, take a resolution
\[
0 \rightarrow E' \rightarrow L'_1 \rightarrow L'_0 \rightarrow P^{\vee}
\rightarrow M^* \rightarrow 0.
\]
Then let $E = E'{}^{\vee}$ and let ${\mathcal E}$ be the associated sheaf.  One
checks easily that ${\mathcal E}$ is extraverti satisfying $T$ and gives rise to
$M$, $P$, $\alpha$ as required.

\bigskip
\noindent
{\bf Corollary 3.3.}  {\em Let $X$ be an arithmetically Gorenstein scheme of
dimension $3$.  Then to each curve $C \subseteq X$ is associated a maximal
Cohen--Macaulay module $P$ over the ring $S = H_*^0({\mathcal O}_X)$, and a map
$\alpha: P^{\vee} \rightarrow M^* \rightarrow 0$, where $M = H_*^1({\mathcal
I}_C)$ is the Rao module of $C$.  Two curves $C_1$ and $C_2$ are in the same
biliaison equivalence class if and only if (up to shift) $M_1 \cong M_2$,
$P_1,P_2$ are stably equivalent, and $\alpha_1,\alpha_2$ are compatible with
this isomorphism on a stable equivalence.  Furthermore, every triple
$(M,P,\alpha)$ as in {\em (3.2c)} occurs for some biliaison class of curves.}

\bigskip
\noindent
{\em Proof.} Given $C$, take any extraverti sheaf ${\mathcal E}$ in the psi
equivalence class of ${\mathcal I}_C$ $(1.12)$.  Then we obtain the associated $P$
and
$\alpha$ by $(3.2)$.  Now if two curves $C_1$ and $C_2$ have compatible $M_i$,
$P_i$, and $\alpha_i$, up to shift, the corresponding sheaves ${\mathcal
E}_1$ and ${\mathcal E}_2$ are stably equivalent, up to shift, by the theorem. 
This in turn implies $C_1$ and $C_2$ are in the same biliaison equivalence class
by $(2.4)$, cf.\ also $(2.14)$.

\bigskip
\noindent
{\bf Remark 3.4.} In the special case $X = {\mathbb P}^3$, the maximal
Cohen--Macaulay module is free, since $S$ is a regular ring.  Thus the condition
of $(3.3)$ boils down to $M_1 \cong M_2$, which is the original theorem of Rao. 
Note that the statement ``every maximal Cohen--Macaulay module is free'' over a
local Cohen--Macaulay ring implies that every module has a finite free
resolution, and so the ring is regular.  Thus ${\mathbb P}^3$ is the only
arithmetically Gorenstein scheme for which the isomorphism of Rao modules is
sufficient to imply biliaison equivalence of curves.

\bigskip
For a two-dimensional arithmetically Gorenstein scheme, there is a much simpler
analogous result.

\bigskip
\noindent
{\bf Proposition 3.5.} {\em Let $X$ be an arithmetically Gorenstein scheme of
dimension $2$.  There is a natural one-to-one correspondence between biliaison
classes of $0$-dimensional closed subschemes $Z$ and stable equivalence classes
of maximal Cohen--Macaulay modules over the ring $S = H_*^0({\mathcal O}_X)$
whose associated sheaf ${\mathcal E}$ is orientable.}

\bigskip
\noindent
{\em  Proof.} To any $0$-dimensional subscheme $Z$ we associate an extraverti
sheaf ${\mathcal E}$ in the psi-equivalence class of ${\mathcal I}_Z$ $(1.11)$. 
Then ${\mathcal E}$ is reflexive $(1.13)$ and hence locally Cohen--Macaulay on
$X$.  Since $X$ is arithmetically Gorenstein, the condition $H_*^1({\mathcal
E}^{\vee}) = 0$ implies by Serre duality that $H_*^1({\mathcal E}) = 0$.  Hence
$E = H_*^0({\mathcal E})$ is a maximal Cohen--Macaulay module over $S$.  Now
$(2.14)$ shows that biliaison equivalence classes of $Z$ correspond to stable
equivalence classes of ${\mathcal E}$ and hence of $E$.  The condition orientable
just requires that ${\mathcal E}$ be locally free in codimension $1$ and its
determinant be isomorphic to ${\mathcal O}_X(\ell)$ for some $\ell$ on $X$ minus
a finite number of points.


\begin{thebibliography}{10}

\bibitem[1]{BBM} Ballico, E., Bolondi, G., and Migliore, J. C., The
Lazarsfeld--Rao problem for liaison classes of two-codimensional subschemes of
${\mathbb P}^n$, {\em Amer. J. Math.} {\bf 113} (1991) 117--128.

\bibitem[2]{BM} Bolondi, G., and Migliore, J. C., The Lazarsfeld--Rao property
on an arithmetically Gorenstein variety, {\em Manusc. Math.} {\bf 78} (1993)
347--368.

\bibitem[3]{GD} Hartshorne, R., Generalized divisors on Gorenstein schemes, {\em
K-theory} {\bf 8} (1994) 287--339.

\bibitem[4]{GDB} Hartshorne, R., Generalized divisors and biliaison, preprint.

\bibitem[5]{HMDP1} Hartshorne, R., Martin--Deschamps, M., and Perrin, D., Un
th\'eor\`eme de Rao pour les familles de courbes gauches, {\em J. Pure Appl.
Algebra} {\bf 155} (2001) 53--76.

\bibitem[6]{HMDP2} Hartshorne, R., Martin--Deschamps, M., and Perrin, D.,
Construction de familles minimales de courbes gauches, {\em Pacific J. Math.}
{\bf 194} (2000) 97--116.

\bibitem[7]{H} Horrocks, G., Vector bundles on the punctured spectrum of a local
ring, {\em Proc. Lond. Math. Soc.} {\bf 14} (1964) 689--713.

\bibitem[8]{LR} Lazarsfeld, R., and Rao, P., Linkage of general curves of large
degree, {\em Springer LNM} {\bf 997} (1983) 267--289.

\bibitem[9]{MDP} Martin--Deschamps, M., and Perrin, D., Sur la classification
des courbes gauches, {\em Ast\'erisque} {\bf 184--185} (1990).

\bibitem[10]{M} Migliore, J. C., {\em Introduction to Liaison Theory and
Deficiency Modules}, Birkh\"auser, Boston (1998).

\bibitem[11]{Na} Nagel, U., Even liaison classes generated by Gorenstein
linkage, {\em J. Algebra} {\bf 209} (1998) 543--584.

\bibitem[12]{No} Nollet, S., Even linkage classes, {\em Trans. AMS} {\bf 348}
(1996) 1137--1162.

\bibitem[13]{R1} Rao, P., Liaison among curves in ${\mathbb P}^3$, {\em Invent.
Math.} {\bf 50} (1979) 205--217.

\bibitem[14]{R2} Rao, P., Liaison equivalence classes, {\em Math. Ann.} {\bf
258} (1981) 169--173.

\bibitem[15]{S} Strano, R., Biliaison classes of curves in ${\mathbb P}^3$.


\end{thebibliography}
\end{document}